\documentclass[a4paper,11pt,leqno]{article}

\usepackage{amsmath,amssymb,amsfonts,latexsym,stmaryrd}

\usepackage{theorem}

\usepackage{natbib}
\usepackage{apalike}

\usepackage[english]{babel}

\usepackage[margin=2.6cm,nohead]{geometry}

\usepackage{enumerate}
\usepackage{multirow}
\usepackage{soul}
\usepackage{bbding}
\usepackage{pifont}
\usepackage{setspace}
\usepackage{comment}
\usepackage{lingmacros}
\usepackage[all]{xy}
\usepackage{bm}
\usepackage{verbatim}
\usepackage[normalem]{ulem}
\usepackage{graphicx,xcolor}
\usepackage{boxedminipage}
\usepackage{fancyhdr}
\usepackage{appendix}
\usepackage{epigraph}
\usepackage{tikz}
\usetikzlibrary{shapes.geometric,decorations.pathreplacing}

\usepackage{subcaption}

\theorembodyfont{\rmfamily}

\newtheorem{Def}{Definition}

\begin{document}
	
	\author{Paul Gorbow and {\O}ystein Linnebo}
	\title{Feasible constructivism} %Feasibilism as a coherent view below intuitionism; Feasibilism as a coherent sub-intuitionism
	\date{Draft of 21 January 2026}
	
	\hyphenation{quanti-fi-ca-tion}

	\maketitle

	\begin{abstract}
		\noindent Dummett's argument for intuitionism is well known. There is a concern that the argument proves too much, specifically, that it supports the extreme and apparently incoherent position of \textit{strict finitism}. The central question is how to explicate the notion that it is \textit{possible in practice} to construct an arithmetical term or verify a  statement. The strict finitist answer is plagued by the sorites paradox. We propose and develop \textit{feasibilism} as a more plausible view, where computational feasibility, as captured by the class of polynomial-time problems, yields a robust and expedient explication of ``possible in practice''. In this approach, the complexity is bounded by a polynomial function of the input size, rather than bounded by a constant (as in strict finitism), thus resolving the sorites issues. We show that a system of strictly bounded arithmetic, introduced by Sam Buss, precisely formalizes the feasibilist view so as to %strongly 
		satisfy Dummett's requirements.
	\end{abstract}
	
	\begin{spacing}{1.5}
	\end{spacing}
		
		\section{Introduction }
		
		Intuitionism is often regarded as a painful restriction of classical mathematics. Hilbert, for example, characterizes intuitionism as ``proscribing [\ldots] to the boxer  the use of his fists'' \cite[p. 476]{Hilbert1927}. 
		Despite this widespread sentiment and Hilbert's scathing remark, we propose to investigate whether there are coherent and interesting views that are yet more restrictive. 
		
		We owe an explanation, clearly, of why sub-intuitionistic views are worthy of  attention. Wasn't it bad enough to deprive the boxer of the use of his fists? We see two reasons to be interested in such views. 
		
		One reason %to be interested in very weak views 
		is encapsulated in the \cite{Feferman1998} observation that  
		``a little bit goes a long way''. It turns out we can do surprisingly much mathematics with comparatively modest resources. Moreover, when we prove a result from weak assumptions, we learn something that a classical proof from stronger assumptions would not have taught us.\footnote{This observation is frequently attributed to Georg Kreisel, see, e.g., \cite{Feferman1996}.}  It can therefore be interesting to ask what we can prove with very weak assumptions. It is important to stress that this question can be interesting even when we are prepared to make stronger assumptions. 
		Curiosity about how well one might live on very modest resources must not be conflated with a desire to live one's entire life that frugally. %It might be a salutary way to focus the mind on the essentials in life. 
		
		A second reason to be interested in sub-intuitionistic views, which will be our primary concern in this paper, is a well-known worry that arguments for intuitionism prove too much. If successful at all, these arguments seem to establish a view more radical than intuitionism. Since sub-intuitionistic views are typically found deeply problematic, this would be a troubling result for proponents of intuitionism. 
		
		The only sub-intuitionistic view that is widely known, at least among philosophers, is \textit{strict finitism}. This is the view that there are only finitely many numbers, such that there an upper bound on how large numbers there are. As we will see shortly, this view is weak to the point of incoherence. 
		
		Our aim is to articulate a clearly coherent and far more interesting sub-intuitionistic view, called  \textit{feasibilism}, which is based on polynomial-time computability. 
		In terms of logical strength, we have: 
		\vspace{2mm}
		\begin{equation*}
			\text{strict finitism} \ < \ \text{feasibilism} \ < \ \text{intuitionism} 
		\end{equation*}
		%{\color{blue}In short, we aim to articulate a radical form of constructivism that is \ldots eh, more feasible than strict finitism. %feasible, both in a technical sense to be explained shortly and in a more everyday sense.
			%}  
		
		Here is our plan. First, we present the concern about proving too much, as it applies to Dummett's arguments.\footnote{A similar concern presumably arises for Brouwer's arguments. We will not pursue that question.} Then, we introduce feasibilism and explain how this view addresses the concern. Finally, we  develop feasibilism properly and show that it meets Dummett's requirements.

		\section{Does Dummett's argument prove too much?}
		
		Dummett’s argument for intuitionistic logic proceeds from a meaning-theoretic constraint, that understanding a statement consists in knowing what would count as a construction and proof of it, to the conclusion that intuitionistic reasoning is justified, but that classical inferences such as the law of excluded middle are not justified in general, e.g. \cite{Dummett1978}. In typical fashion, though, \cite{Dummett1975} articulates a concern about this argument for intuitionism. The concern has the following structure: 
		\begin{enumerate}[(1)]
			\item Dummett’s argument for intuitionism supports strict finitism as well. 
			
			\item 	Strict finitism is incoherent. 
			
			\item Therefore, we have identified a reason to reject Dummett’s argument. 
		\end{enumerate}
		Let us spell out the two premises, (1) and (2), in that order. 
		
		%Beginning with (1), 
		According to Dummett, % characterizes a constructivist view of mathematics as follows
		\begin{quote}
			[c]onstructivist philosophies of mathematics insist that the meanings of all terms [\ldots] %, including logical constants, 
			appearing in mathematical statements must be given in relation to \emph{constructions which we are capable of effecting}, and of \emph{our capacity to recognise} such constructions as providing proofs of those  statements \cite[p.~301, our italics]{Dummett1975}
		\end{quote}
		Thus, Dummett's constructivism centers on requirement that an account of meaning be given in terms of what we are ``capable of effecting'' or recognizing.
		\begin{comment}
			\begin{quote}
				Constructivist philosophies of mathematics insist that the meanings of all terms, including logical constants, appearing in mathematical statements must be given in relation to \textit{constructions which we are capable of effecting}, and of \textit{our capacity to recognise such constructions} as providing proofs of those statements; and, further, that the principles of reasoning which, in assessing the cogency of such proofs, we acknowledge as valid must be justifiable in terms of the meanings of the logical constants and
				of other expressions as so given.    (our italics) 
			\end{quote}
		\end{comment}
		He defends this constraint by reflecting on how we can learn the language of mathematics:
		\begin{quote}
			We learn, and can only learn, their meanings by a training in their use; and that means a training in \textit{effecting mathematical constructions}, and in \textit{recording them} within the language of mathematics. There is no means by which we could derive from such a training a grasp of anything transcending it, such as a notion of truth and falsity for mathematical statements independent of \textit{our means of recognising} their truth-values.    (our italics) \cite[p. 301]{Dummett1975}
		\end{quote}
		
		Assume we grant Dummett's call for a  constructivist account of meaning. As he observes, this will have consequences for logic, because
		\begin{quote}
			%    Constructivist philosophies of mathematics insist that the meanings of all terms, including logical constants, appearing in mathematical statements must be given in relation to \textit{constructions which we are capable of effecting}, and of \textit{our capacity to recognise such constructions} as providing proofs of those statements; and, further, that 
			the principles of reasoning which [\ldots] %, in assessing the cogency of such proofs, 
			we acknowledge as valid must be justifiable in terms of the meanings of the logical constants and of other expressions as so given.    
		\end{quote}
		After all, Dummett's constructivism %t ``job description'' 
		entails that any theorem of the form $\exists x \varphi(x)$ must be backed up by a \emph{manifestation}, that is, an effective construction of an object $a$ together with a verification that $\varphi(a)$. Moreover, it entails that only intuitionistic reasoning is justified. %We take the view thus expressed as  
		In the case of arithmetic, \cite[p. 36]{Dummett1977} takes this view to be formalized by Heyting Arithmetic ($\mathsf{HA}$). This is quite natural, since $\mathsf{HA}$ has the well-known existence property.
		
		So far, so good---Dummett thinks. Now for his concern. Crucial to Dummett's argument are the notions of ``constructions we are capable of effecting'' and ``our capacity to recognize''. How should these notions be understood? Idealizing radically, intuitionism understands both notions in terms of what we can do \textit{in principle}. Is this idealization warranted? Why not restrict to what we can do \textit{in practice}? That is, why not insist 
		\begin{quote}
			that the meanings of our terms must be given by reference to \emph{constructions which we can in practice carry out}, and to criteria of correct proof on which \emph{we are in practice prepared to rely} (ibid., our italics)
		\end{quote}
		%We will be interested in an even stricter view: the only math object that exist are those that we can in practice construct; the only truths, those that we can in practice prove.
		This interpretation of the Dummett's constraint on an of meaning calls for a view stricter than intuitionism. And that is dangerous for Dummett. For the usual explication of the more radical view is \textit{strict finitism}. 
		
		As mentioned, if strict finitism holds there is an upper bound on how large numbers there are. As \cite{Dummett1975} puts it, ``by `natural number' must be understood a number which we are in practice capable of representing'' (p. 303). Thus, strict finitists often deny that there is a number such as $10^{10^{10}}$, on the grounds that we cannot in fact produce the associated canonical numeral.
		Clearly, this is an extremely radical view---so radical, we claim, that it is problematic. 
		
		Are the natural numbers closed under the usual arithmetical operations, including that of \textit{successor}? Strict finitists face a dilemma. Suppose they deny that closure property. That would cripple mathematics---a fate even worse than being deprived of the use of one's fists.
		
		Alternatively, suppose that strict finitists grant the totality of the successor function. Then they would need to reconcile that admission with their view that only small natural numbers exist. As \cite{Dummett1975} observes, this leads to a version of the sorites paradox that he calls Wang's paradox: 
		\begin{itemize}
			\item 0 is small
			\item if $n$ is small, then $n'$ is small
		\end{itemize}
		Dummett's diagnosis is that ``$n$ is small'', regardless if interpreted as ``$n$ grains of sand
		are too few to make a heap'' or as ``it is possible in practice to write down the Arabic numeral for $n$'', is vague. Based on this vagueness, \cite{Dummett1975} gives a detailed argument that any formalization of strict finitism is incoherent. On the other hand, \cite{Wright1982} claims to overcome this hurdle. Be that as it may. We believe it is already a serious awkwardness of strict finitism that it needs to distinguish between numbers that are and are not ``small''. In addition to Wang's paradox, what counts as small will depend on parochial facts %is subject to contingencies 
		concerning what kinds of creatures we are and the kind of environment that we inhabit. %]} makes a cut-off in its formalization of the vague notion of constructibility and verifiability in practice.
	
	Thus, on either horn of the dilemma, strict finitism appears deeply problematic.

	\section{Feasibilism as a coherent and interesting view} %more plausible alternative
	
	There is, however, a coherent and interesting sub-intuitionistic view---which is also a more plausible end-point of Dummett's arguments than strict finitism. 
	
	Our proposal, which we call \emph{feasibilism}, formalizes constructibility and verifiability in practice as constructible and verifiable by a polynomial time algorithm with respect to the complexity of the dependent variables. In particular, this entails that for atomic $\varphi$, $\forall x \exists y \varphi(x, y)$ is verified if, and only if, there is a polynomial-time Turing machine, such that for every input term $s$ it outputs a term $t$ and an intuitionistc proof of $\varphi(s, t)$. The identification of feasible computation with polynomial time computation is deeply embedded as a methodological principle in the field of computer science. It goes by the name \emph{Cobham's thesis}, and is introduced and defended in \cite{Cobham1965}. 
	
	We claim that feasibilism allays  Dummett's concern. For we appear to have: 
	\begin{enumerate}[(1$'$)]
		\item Dummett’s argument for intuitionism supports feasibilism, not strict finitism.   
		
		\item 	Feasibilism is coherent.
		
		\item Therefore, we have not identified a reason to reject Dummett’s argument. 
	\end{enumerate}
	
	Whether we should accept Dummett’s argument is of course a  different question. For the record: we are skeptical. Though that will not play any role in this note.
	
	Let us defend the two premises of the argument. %, (1$'$) and (2$'$). 
	We begin with the coherence of feasibilism, (2$'$), as this is most straightforward. 
	Feasibilism clearly avoids the problems of strict finitism. It grants the totality of the successor function. In this way, it also avoids the sorites paradox. As further evidence for its coherence---and substantial theoretical interest---we spell out the view more properly in the next section. 
	
	%We do claim, however, that feasibilism is interesting because it is a coherent and conceptually fairly robust view, which makes only the weakest kind of assumptions. 
	
	Turning to the first premise, (1$'$), the crux is how to explicate ``which we are capable of effecting''. \textit{What is the right level of idealization}? Here are the three main options and the views to which they give rise. 
	\begin{itemize}
		\item none: an absolutely resource-bounded computer \quad $\rightsquigarrow$ strict finitism
		
		% specific in an unscientific way: why human rather than Martian? which humans? why not human-aided-by-supercomputer? 
		
		\item a relatively resource-bounded computer \qquad $\rightsquigarrow$ feasibilism
		
		% This matches the letter of Dummett's more demanding formulation. 
		
		\item a resource-unbounded computer \qquad $\rightsquigarrow$ intuitionism
	\end{itemize}
	%\begin{itemize}
	%    \item none: an actual human mathematician \quad $\rightsquigarrow$ strict finitism\vspace{1mm}
	%    \item or aided by a ``realistic'' computer \qquad $\rightsquigarrow$ feasibilism\vspace{1mm}
	%    \item or aided by an idealized computer \qquad $\rightsquigarrow$ intuitionism
	%\end{itemize}

	We find our proposed view of feasibilism to be precisely what is supported by Dummett’s
	argument on learning meanings by training in their use. Naturally, the complexity of the training process depends on the complexity of the term whose meaning is to be learned. As feasibilism formalizes this very dependence in accordance with the canonical complexity class of feasibility, namely polynomial time, this view is theoretically robust and dissolves the discord between vagueness and arbitrary cut-offs.
	
	As noted, Dummett’s characterization of constructivism hinges on the phrase ``constructions which we are capable of effecting'', and the crucial issue is how that capacity is to be explicated. One option is to interpret ``capable'' in the most restrictive possible sense, identifying it with the resources of a fixed, absolutely bounded computational agent. This yields strict finitism: only constructions that fit under a hard, non-negotiable resource ceiling count as meaningful. But as Dummett himself stresses, this approach inevitably reintroduces the vagueness of Wang’s paradox at the formal level. Moreover, any such hard resource ceiling must be fixed independently of the mathematical tasks under consideration. As a result, it does not track any general structural feature of mathematical practice or reasoning, but merely reflects the limitations of a particular agent or context. From a scientific point of view, this is problematic: science typically seeks principles with general explanatory and predictive power, rather than theories whose validity depends on the contingent capacities of a specific finite system. 
	
	At the opposite extreme lies an interpretation in which the relevant agent is resource-unbounded. Here, every construction permitted by an idealized, non-finitary computing power becomes legitimate---exactly the stance captured by intuitionistic arithmetic (HA). While this dissolves the finitist worries, it arguably goes well beyond Dummett’s more demanding semantic requirements. The emphasis on training and recognition suggests something more modest and operational that can reasonably be carried out in practice.
	
	Feasibilism occupies the principled middle ground. If ``capable of effecting'' is understood relative to a \emph{resource-bounded} but \emph{unboundedly extensible} computational model, then polynomial-time computation emerges as the natural explication. It avoids the arbitrariness of strict finitism while respecting the semantic constraints highlighted by Dummett. Crucially, this grounded middle-way preserves both the spirit and the letter of Dummett's semantic constraint: constructions must be possible \emph{in practice}, but ``in practice'' scales with the complexity of the input in a rule-governed, non-arbitrary manner. This relative interpretation of practical possibility frees the notion from dependence on any special resource constraint conditioned on a particular agent. Thus we obtain a robust notion that is methodologically expedient to the scientific aim of general explanatory power. In the next section we proceed to fully explicate feasibilism as a robust scientific view, by providing a formal framework which strongly satisfies Dummett's demands and articulates them with precision.

	\section{Feasibilism explicated and %gainfully 
		put to use}\label{sec:feasibilism-framework}
	
	As just observed, feasibilism avoids Wang's paradox, which is a first step towards establishing its coherence. But more work remains. We need to explicate feasibilism precisely along the motivating idea that only what is polynomial-time computable is feasible. And we need to show that feasibilism satisfies Dummett's constraint on an account of meaning, such that this approach can be used as a radical form of constructivism. We rely heavily on \cite{Buss1986}, a systematic logical characterization of polynomial time complexity, which will be the default reference throughout this section.\footnote{Since Buss's seminal work, there have been significant developments on characterizing computational complexity classes in logical terms, particularly for the polynomial-time class. Still, Buss's characterization stands out as the most straightforward one for philosophical purposes, because it is a relatively simple classical first-order fragment of the familiar Peano Arithmetic, $\mathsf{PA}$. It was followed by a second-order characterization developed in \cite{Leivant1991}. Recent developments have focused on characterizing polynomial-time complexity in \textit{linear logic}, a sub-intuitionistic logic that restricts the ability to re-use assumptions. In complexity-theoretic terms this is natural, because copying consumes real computational resources. Notable examples of that line of research are \cite{CrossleySeely1994}, \cite{BrunelTerui2010}, and \cite{DalLagoHofmann2011}.} 
	
	First, we articulate the feasibilist view precisely, and axiomatize a first-order theory $\mathsf{FA}$ of Feasible Arithmetic (\S \ref{sec:FA}) intended to capture this view. Second, we explain how this system naturally characterizes polynomial-time computability, our core formalization of feasibility (\S \ref{sec:total_polytime}).  Third, we give a realizability semantics for feasibilism within $\mathsf{FA}$ (\S \ref{sec:realizability}), which shows the precise and strong sense in which $\mathsf{FA}$ satisfies Dummett's semantic constraints. In particular, the relevant semantic criteria of ``construction we are capable of effecting'' and of ``our capacity to recognize such constructions as providing proofs'', as well as the requirement that ``the principles of reasoning which \dots we acknowledge as valid must be justifiable in terms of the meanings'' are made precise within $\mathsf{FA}$ in accordance with Cobham’s thesis.

	\subsection{The system of feasible arithmetic}\label{sec:FA}
	
	We take \emph{feasibilism} to be the view that an arithmetical construction is valid precisely when it can be carried out in \emph{polynomial time} in the size of its input. Concretely, given a term $t(\vec{x})$, a computation of the value of $t$ relative to particular arguments $\vec{a}$ is required to run in time bounded by a polynomial in the total binary length $|\vec{a}|$ of $\vec{a}$. More generally, in order to address Dummett’s requirements, we extend this notion from mere term evaluation to the construction of \emph{realizers} for arithmetical statements. Thus, for a formula $\varphi(\vec{x})$, feasibilism demands that, for each input $\vec{a}$, $\varphi(\vec{a})$ holds if, and only if, there exists a realizer $r$ for $\varphi(\vec{a})$ that can be computed in polynomial time in $|\vec{a}|$; and that there is likewise a polynomial-time procedure which decides whether such an $r$ is indeed a realizer for $\varphi(\vec{a})$. In short, both the construction of mathematical meaning and the recognition of correctness are taken to be feasible exactly when they are achievable within polynomial-time complexity. Moreover, feasibilism demands that any valid reasoning must preserve these conditions on meaning and recognition.
	
	We find Buss's first-order classical system $S^1_2$ to be a system that has precisely enough strength to meet these requirements, and we suggest to name it $\mathsf{FA}$, for ``Feasible Arithmetic''. The language and axioms of $\mathsf{FA}$ are designed to allow the system to express precisely bounded quantifiers that suffice to logically demarcate certain relevant complexity classes. Moreover, it has a restricted induction schema that is precisely strong enough to prove the recursion theorem for polynomial-time recursion. From this basis, Buss is able to give a natural proof-theoretic characterization of the class of polynomial-time computable functions. As the notion of polynomial-time computable function is complexity-theoretic, and not directly expressible in the language of first-order logic, Buss's logical characterization carries significant depth. This bridge between complexity-theory and logic is precisely what we need for a system of feasible arithmetic.
	
	The first-order language of $\mathsf{FA}$ includes all standard logical symbols
	$\land,\,\lor,\,\neg,\,\supset,\,\exists,\,\forall$ together with parentheses,
	and the nonlogical symbols
	\[
	S,\,0,\,+,\,\cdot,\,|x|,\,\bigl\lfloor\tfrac{x}{2}\bigr\rfloor,\,\#,\text{ and } \leq.
	\]
	These nonlogical symbols are interpreted over the natural numbers.  Here $S$ denotes the successor function,
	$0$ the zero constant, $+$ and $\cdot$ the usual addition and multiplication,
	and $\leq$ the standard order relation.  The term $|x|$ designates the
	\emph{length} of the binary representation of $x$, that is,
	\[
	|x|=\lfloor \log_2(x+1)\rfloor.
	\]
	The term $\lfloor\tfrac{x}{2}\rfloor$ designates the greatest integer less than or equal
	to $x/2$. This operation right-shifts the binary representation of $x$. For example, the number $10011$ is mapped to $1001$, and $1010$ is mapped to $101$. Finally, the term $x\# y$, pronounced ``$x$ smash $y$'', designates
	\[
	x\# y = 2^{\,|x|\cdot |y|}.
	\]
	This operation turns out to be methodologically expedient, because it allows one to express that a quantifier is bounded by $2^{p(|x|)}$, where $p$ is a polynomial. In the complexity-theoretic correspondence, this bounds the size $|y|$ of the quantified variable $y$ by a polynomial of the size $|x|$ of the input $x$, ensuring that the binary representation of $y$ can be read in polynomial time in $|x|$.
	
	Let $\Delta^b_0=\Sigma^b_0=\Pi^b_0$ denote the class of
	\emph{sharply bounded} formulas, i.e.\ formulas whose quantifiers
	are all of the form $(\forall x\le |t|)\,\varphi$ or
	$(\exists x\le |t|)\,\varphi$.
	The next levels $\Sigma^b_1$ and $\Pi^b_1$ are obtained by adding
	one bounded alternation of quantifiers.
	
	\begin{Def}
		\begin{enumerate}
			\item $\Sigma^b_1$ is the smallest set of formulas satisfying:
			\begin{enumerate}
				\item $\Sigma^b_1\supseteq \Delta^b_0$;
				\item if $\varphi\in\Sigma^b_1$, then $(\exists x\le t)\,\varphi$
				and $(\forall x\le |t|)\,\varphi$ are in $\Sigma^b_1$;
				\item if $\varphi,\psi\in\Sigma^b_1$, then $\varphi\wedge \psi$ and $\varphi\vee \psi$
				are in $\Sigma^b_1$;
				\item if $\varphi\in\Sigma^b_1$ and $\psi\in\Pi^b_1$, then $\neg \psi$
				and $\psi\supset \varphi$ are in $\Sigma^b_1$.
			\end{enumerate}			
			\item $\Pi^b_1$ is defined dually as the smallest set of formulas satisfying:
			\begin{enumerate}
				\item $\Pi^b_1\supseteq \Delta^b_0$;
				\item if $\varphi\in\Pi^b_1$, then $(\forall x\le t)\,\varphi$
				and $(\exists x\le |t|)\,\varphi$ are in $\Pi^b_1$;
				\item if $\varphi,\psi\in\Pi^b_1$, then $\varphi\wedge \psi$ and $\varphi\vee \psi$
				are in $\Pi^b_1$;
				\item if $\varphi\in\Pi^b_1$ and $\psi\in\Sigma^b_1$, then $\neg \psi$
				and $\psi\supset \varphi$ are in $\Pi^b_1$.
			\end{enumerate}
		\end{enumerate}
	\end{Def}
	
	This can be extended to a full-blown hierarchy, defining $\Sigma^b_n$ and $\Pi^b_n$, for all $n \in \mathbb{N}$. The union of that whole hierarchy is equivalent to the class $\Delta_0$ of formulas with only bounded quantifiers, which forms the base of the well-known {\em arithmetical hierarchy}. Thus, we are here dealing with a much more restricted hierarchy, where sharply bounded and bounded quantifiers play roles analogous to bounded and unbounded quantifiers, respectively, in the arithmetical hierarchy. Interestingly, $\Delta^b_0$ is contained in the class of polynomial-time decidable predicates, and $\Sigma^b_1$ is precisely the class of non-deterministic polynomial-time decidable predicates.\footnote{The latter is a non-trivial result. See Buss's \S 1, and its Theorem 8 in particular, for formal details. The class of non-deterministic polynomial-time decidable predicates can be roughly defined as the closure of the class of polynomial-time decidable predicates under the quantifier $\exists y \leq 2^{p(|\vec{x}|)}$, for any freely chosen polynomial $p$. Moreover, the problem of deciding $\exists y \leq 2^{p(|\vec{x}|)} \, R(y, \vec{x})$ can be characterized as the problem of verifying $R(b, \vec{x})$, for a given $b$. Thus, the class of non-deterministic polynomial-time decidable predicates is intuitively the class of polynomial-time verifiable predicates.}
	
	We write $\Delta^b_1$ for formulas equivalent (in
	$\mathsf{FA}$ as axiomatized below) to both a $\Sigma^b_1$ and a $\Pi^b_1$ formula.
	
	Let $\mathsf{FA}$ (\emph{Feasibilist Arithmetic}) be the first-order theory $S^1_2$ introduced in Buss's \S 2.2. The axiomatization consists of the natural inductive equations formalizing the function symbols, and the induction axiom schema $\Sigma^b_1$-$\mathsf{PIND}$: 
	\begin{align*}
		(\varphi(0)\wedge \forall x(\varphi(\lfloor x/2\rfloor)\to \varphi(x)))\to
		\forall x\,\varphi(x) \textnormal{, for all $\varphi\in\Sigma^b_1$.} \tag{$\Sigma^b_1$-$\mathsf{PIND}$}
	\end{align*}
	We may think of this axioms schema as expressing tree-induction on the binary representations of natural numbers. It says that for any $\Sigma^b_1$-formula $\varphi(x)$, if 
	\begin{itemize}
		\item $0$ satisfies $\varphi$, and 
		\item for any number $x$, binary-represented by $b_1 \cdots b_n$, if the number binary-represented by the right-shifted sequence $b_1 \cdots b_{n-1}$ satisfies $\varphi$ then $x$ satisfies $\varphi$,
	\end{itemize}
	then all numbers satisfy $\varphi$.
	
	This weakens the induction schema of $\mathsf{PA}$ in two ways: First, it is severely restricted to the narrow class of $\Sigma^b_1$-formulas. 
	Second, the tree-induction formulation is weaker than the usual formulation in terms of the successor function. Indeed, it is typically harder to prove the conditional $\forall x(\varphi(\lfloor x/2\rfloor)\to \varphi(x))$, than to prove the conditional $\forall x(\varphi(x)\to \varphi(S(x)))$.
	
	As shown in Buss's \S 2.4--5, $\mathsf{FA}$ is strong enough to represent many useful relations and functions. In particular, it defines tuples and a projection function $\beta$, such that $\mathsf{FA}$ proves ``$\beta(i, r) =$ the $i$-th component of the tuple encoded by $r$''; this is explicitly used to define realizability in \S \ref{sec:realizability}.
	
	$\mathsf{FA}$ is a system for reasoning about $\Sigma_1^b$-formulas. Indeed, the induction schema only applies to such formulas, and the results proved about $\mathsf{FA}$ are only concerned with this fragment. Thus, any theorems of higher complexity proved by $\mathsf{FA}$, such as instances of excluded middle, are considered out of scope with regard to our formalization of Feasible Arithmetic.
	
	We have now reached a point where \emph{Feasible Arithmetic} stands fully and
	transparently formalized. The system is presented here with its intended feasibilist motivation: the language isolates exactly the sharply bounded and
	bounded quantifiers needed to capture polynomial-time computation, and the sole
	induction schema $\Sigma^{b}_{1}$\text{-PIND} matches the combinatorial structure of
	polynomial-time recursion. All function symbols are given by
	elementary defining equations, and the induction principle is a
	philosophically natural expression of resource-bounded induction on the relevant computational structure (namely binary-tree representations). It is a concrete fragment of $\mathsf{PA}$, calibrated precisely to
	polynomial-time reasoning.
	
	Moreover, $\mathsf{FA}$ sits inside the intuitionistic system of Heyting Arithmetic, $\mathsf{HA}$, in a technically exact way. Let
	$\mathsf{FA}^{\mathrm{int}}$ be the intuitionistic version of $\mathsf{FA}$---that is, the same
	language and axioms,
	but formulated in intuitionistic logic. Then the classical/intuitionistic gap disappears
	completely at the level of feasible content. By a result of Avigad
	(Theorem~3.17 in~\cite{Avigad2000}), $\mathsf{FA}$ is \emph{conservative} over its
	intuitionistic counterpart for all $\Sigma^{b}_{1}$–formulas. Since the results we rely on in this paper are all formulated only for $\Sigma^{b}_{1}$–formulas, they all hold for $\mathsf{FA}^{\mathrm{int}}$ as well. Moreover, as our formalization of Feasible Arithmetic is only concerned with $\Sigma^{b}_{1}$–theorems, it is, in a
	rigorous sense, a subsystem of $\mathsf{HA}$. This picture is reminiscent of Brouwer’s well-known position that
	classical logic is the logic of the finite \cite[p. 336]{vanHeijenoort1967}. More generally, once one works inside finitely bounded domains (here formalized by bounded quantifiers), classical reasoning ceases to outrun constructive justification.
	
	At this point, a reader might wonder whether the use of strictly bounded quantifiers in $\mathsf{FA}$ reflects an ineliminable dependence on syntax, or merely a natural way to control complexity. Recent work by Tabatabai and Thapen (not yet a published article, we refer to their long conference abstract: \cite{TabatabaiThapen2024}) suggests that the latter is the case. They introduce a weak predicative set theory, PST, with sorts for sets and classes, which is naturally axiomatized by extensionality, empty set, pairing, union, $\Delta_0$-separation for sets, and $\Sigma_1$-comprehension for classes. PST admits an interpretation of $\mathsf{FA}$ on binary strings, while conversely the $\Sigma_1$-fragment of PST restricted to binary strings can be interpreted inside $\mathsf{FA}$. As a result, the provably total $\Sigma_1$-definable functions of PST over binary strings are exactly the polynomial-time computable ones. This shows that feasible computation can be captured structurally---via predicative separation and comprehension---in standard set-theoretic syntax. From the perspective of feasibilism, this result is reassuring: it indicates that the appeal to strict bounds in $\mathsf{FA}$ is only one natural way of enforcing feasibility within a first-order arithmetical setting. Alternative formalisms equivalently achieve the same constraint, reinforcing the case that Feasible Arithmetic is theoretically robust and invariant under syntactic contingencies.

	\subsection{Feasibly computable functions and relations}\label{sec:total_polytime}
	
	A necessary condition on any purported explication of feasibilism is that it proves the totality of precisely the functions that are in fact polynomial-time computable. This condition is indeed satisfied:
	
	By Corollary~8 of Buss's~\S 5.3, the provably total $\Sigma^b_1$-definable functions of
	$\mathsf{FA}$ are exactly the \emph{polynomial-time computable functions}.
	Formally, if $\varphi\in\Sigma^b_1$ and
	\[
	\mathsf{FA}\ \vdash\ \forall x\,\exists!y\, \varphi(x,y),
	\]
	then there exists a polynomial-time computable function $f$ such that
	\[
	\mathbb N\models\forall x\,\varphi(x,f(x)).
	\]
	Conversely, every polynomial-time function is $\Sigma^b_1$-definable in
	$\mathsf{FA}$. 
	
	This theorem shows that $\mathsf{FA}$ formalizes the feasibilist view of functions. Indeed, given any arithmetic function $f(\vec x)$, we have:
	\[
	\mathsf{FA} \vdash \text{``$f$ is }\Sigma^b_1\text{-definable and total''} \quad\Longleftrightarrow\quad
	f(\vec x)\ \text{is polynomial-time computable.}
	\]
	
	Finally, Theorem~9 in Buss’s \S 5.3 shows that the correspondence between
	provability and feasibility extends from functions to relations.
	Given any arithmetic relation $R(\vec x)$, we have:
	\[
	\mathsf{FA}\ \vdash\ \text{``$R$ is }\Delta^b_1\text{-definable''} \quad\Longleftrightarrow\quad
	R(\vec x)\ \text{is polynomial-time decidable.}
	\]
	This is simply what we get when we represent $R$ as a function into \{0,1\}. Thus, $\mathsf{FA}$ also captures exactly the class of
	feasibly computable relations.

	\subsection{Realizing $\Sigma^b_1$-truths}\label{sec:realizability}
	
	Following Buss's \S 5.1, for each $\Sigma^b_1$-formula $\varphi$ we define a predicate
	$\mathsf{Rlz}_\varphi(r,\vec{x})$ expressing that the object~$r$
	\emph{realizes} the $\Sigma^b_1$-formula $\varphi(\vec{x})$.\footnote{Buss calls this object a witness, and denotes this predicate $\mathit{Witness}_\varphi(w, \vec{x})$.} Intuitively, $r$ encodes a construction verifying~$\varphi$. As we shall see, a realizer for $\varphi(\vec{a})$ is computable in polynomial time in the size $|\vec{a}|$ of $\vec{a}$, and the predicate $\mathrm{Rlz}_\varphi(r, \vec{a})$ is likewise decidable in polynomial time in $|\vec{a}|$. These two facts provide an explicit feasibilist satisfaction of Dummett's semantic constraint, which we here crystallize as follows:
	\begin{description}
		\item[Feasible construction:] The meaning of any mathematical formula must be given compositionally by a feasible construction, in relation to the input given for the free variables.
		\item[Feasible verification:] There must be a feasible procedure for recognizing whether any such construction is a proof of the statement, for a given input.
		\item[Valid inference:] Moreover, the rules of inference must feasibly preserve the above conditions.
	\end{description}

	We now proceed to define Buss's realizability semantics. Let $\varphi$ be a $\Sigma^b_1$-formula with free variables among~$\vec{x}$, and let $t$ be a term. $\mathsf{Rlz}_\varphi(r,\vec{x})$ is defined by the following recursion on the logical complexity of~$\varphi$, which can be carried out within $\mathsf{FA}$ (recall that $\beta(i, r) =$ ``the $i$-th component of the tuple encoded by $r$''):
	\begin{enumerate}
		\item If $\varphi$ is $\Delta_0^b$, then
		\[
		\mathsf{Rlz}_\varphi(r,\vec{x})\ \leftrightarrow\ \varphi(\vec{x}).
		\]
		\item If $\varphi \equiv \psi_1 \wedge \psi_2$, then
		\[
		\mathsf{Rlz}_\varphi(r,\vec{x})\ \leftrightarrow\
		\mathsf{Rlz}_{\psi_1}(\beta(1,r),\vec{x})\wedge
		\mathsf{Rlz}_{\psi_2}(\beta(2,r),\vec{x}).
		\]
		
		\item If $\varphi \equiv \psi_1 \vee \psi_2$, then
		\[
		\mathsf{Rlz}_\varphi(r,\vec{x})\ \leftrightarrow\
		\mathsf{Rlz}_{\psi_1}(\beta(1,r),\vec{x})\vee
		\mathsf{Rlz}_{\psi_2}(\beta(2,r),\vec{x}).
		\]
		\item If $\varphi \equiv (\forall y\le |t(\vec{x})|)\,\psi(\vec{x},y)$, then
		\[
		\mathsf{Rlz}_\varphi(r,\vec{x})\ \leftrightarrow\
		\mathsf{Seq}(r)\wedge \mathsf{Len}(r)=|t(\vec{x})|+1\wedge
		\forall y\le |t(\vec{x})|\,\mathsf{Rlz}_\psi(\beta(y+1,r),\vec{x},y),
		\]
		where $\mathsf{Seq}(r)$ and $\mathsf{Len}(r)$ express that $r$
		is a sequence $\langle r_0,\ldots,r_{|t|}\rangle$ of the appropriate length.
		\item If $\varphi \equiv (\exists y\le t(\vec{x}))\,\psi(\vec{x},y)$, then
		\[
		\mathsf{Rlz}_\varphi(r,\vec{x})\ \leftrightarrow\
		\mathsf{Seq}(r)\wedge \mathsf{Len}(r)=2\wedge
		\beta(1,r)\le t(\vec{x})\wedge
		\mathsf{Rlz}_\psi(\beta(2,r),\vec{x},\beta(1,r)).
		\]
		\item If $\varphi\equiv\neg\psi$, then $\mathsf{Rlz}_\varphi$ is defined by the usual
		syntactic transformation of~$\psi$ pushing the negation down the structure to obtain an equivalent form handled by cases~(1)--(5). For example, $\neg (\psi_1 \wedge \psi_2)$ is transformed to $\neg \psi_1 \vee \neg \psi_2$, which is handled by case (3).\footnote{Thus, negation is ostensibly handled differently than in ordinary realizability, where it is typically treated as a special case of implication, and where an implication is realized by a function that maps a realizer of the antecedent to a realizer of the consequent. However, as shown by Buss's Theorem 5 in his \S 5.2, his realizability semantics implicitly provides such a function for any provable implication; and moreover, it is polynomial-time computable. }
	\end{enumerate}
	
	Thus, $r$ encodes a finite sequence of partial realizers and their verifications, corresponding to the bounded quantifier structure
	of~$\varphi$.

	A central feature of $\mathsf{FA}$ is that
	$\Sigma^b_1$-statements are provably equivalent to the
	existence of \emph{feasible realizers}.  Buss’s Proposition~3 in his \S 5.1
	shows that for any $\Sigma^b_1$ formula $\varphi(\vec{x})$
	\[
	\varphi(x)\;\equiv\;
	(\exists y\le t(\vec{x}))\,R(\vec{x},y),
	\]
	where $R\in\Delta^b_0$ and $t$ is a term, we have:
	\begin{equation*}\label{eq:witnessing}
		\mathsf{FA}\ \vdash\
		\varphi(\vec{x})\ \leftrightarrow\ 
		\exists r\;\mathsf{Rlz}_\varphi(r, \vec{x}).
	\end{equation*}
	Moreover, there is a polynomial-time function $g_\varphi(\vec{x})$ (definable in $\mathsf{FA}$), such that 
	\[\mathsf{FA} \vdash \varphi(\vec{x}) \leftrightarrow \mathsf{Rlz}_\varphi(g_\varphi(\vec{x}), \vec{x}).\]
	
	Thus, every provable $\Sigma^b_1$-truth is realized,
	and conversely the existence of a realizer ensures the truth of the formula. Note that we have now satisfied the constraint of feasible construction, because the meaning of $\varphi$ has been given within $\mathsf{FA}$ in terms of necessary and sufficient compositional conditions (expressed in the definition of $\mathsf{Rlz}_\varphi(r, \vec{x})$) on a construction $r$, which is feasible.
	
	Moreover, the demand for feasible verification is satisfied by Buss’s Proposition~4 in his \S 5.1, which asserts that for each $\Sigma^b_1$-formula $\varphi$,
	the predicate $\mathsf{Rlz}_\varphi(r, \vec{x})$ is itself \emph{polynomial-time
		decidable}.
	
	Hence, every $\Sigma^b_1$-formula provable in $\mathsf{FA}$ is
	accompanied not only by a feasibly computable construction~$r$,
	but also by a feasible \emph{verification} procedure verifying
	that $r$ is indeed a realizer:  
	\[
	\text{provability}\quad\Longleftrightarrow\quad
	\text{feasible construction + feasible verification}.
	\]
	
	Lastly, it follows from Buss's Theorem 5 in his \S 5.2 that inference in $\mathsf{FA}$ feasibly preserves the constructions and verifications of this realizability semantics, thus satisfying the demand for valid inference. More precisely, this theorem entails that if $\varphi(\vec{x})$ is provable in $\mathsf{FA}$ from the assumption $\psi(\vec{x})$, then there is a polynomial-time function $f(y, \vec{x})$ (definable in $\mathsf{FA}$), such that for any realizer $r$ of $\psi(\vec{a})$, we have that $f(r, \vec{a})$ realizes $\varphi(\vec{a})$.
	
	Thus, this framework complies remarkably well with Dummett's semantic constraint. 
	The meaningful mathematical statements are characterized as the $\Sigma^b_1$-statements. Moreover, any relevant construction and accompanying recognition for a formula $\varphi$ are precisely represented within $\mathsf{FA}$ by a realizer and by the predicate $\mathsf{Rlz}_\varphi$, respectively; and inference in $\mathsf{FA}$ is tracked by polynomial-time operations on such realizers. This completes the task of showing that $\mathsf{FA}$ provides a natural, precise, and highly expedient system for explicating feasibilism as a constructivist view aligned with Cobham's thesis.

\end{document}